\documentclass[11pt,reqno]{amsart}

\usepackage{lineno,hyperref}

\makeatletter
\newcommand{\inlineitem}[1][]{%
	\ifnum\enit@type=\tw@
	{\descriptionlabel{#1}}
	\hspace{\labelsep}%
	\else
	\ifnum\enit@type=\z@
	\refstepcounter{\@listctr}\fi
	\quad\@itemlabel\hspace{\labelsep}%
	\fi}
\makeatother

\makeatletter
\g@addto@macro{\endabstract}{\@setabstract}
\newcommand{\authorfootnotes}{\renewcommand\thefootnote{\@fnsymbol\c@footnote}}%
\makeatother

\usepackage{graphicx,epstopdf}
\usepackage{amsfonts,epsfig,fancyhdr,graphics, hyperref,amsmath,amssymb}
\usepackage[inline]{enumitem}
\usepackage{mathtools}
\usepackage{tikz}
\usepackage{amsthm}
\usepackage{cleveref}
\usepackage[colorinlistoftodos]{todonotes}

\theoremstyle{plain}
\newtheorem{theorem}{Theorem}[section]
\newtheorem{corollary}[theorem]{Corollary}
\newtheorem{proposition}[theorem]{Proposition}

\theoremstyle{definition}
\newtheorem{definition}[theorem]{Definition}
\newtheorem{remark}[theorem]{Remark}
\newtheorem{example}[theorem]{Example}

\newcommand{\cT}{\mathcal T}

\newcommand{\dist}{\mathrm{dist}}

\begin{document}
	
	\begin{center}
		\LARGE 
		Kemeny's constant and Wiener index on trees \par \bigskip
		
		\normalsize
		\authorfootnotes
		Jihyeug Jang\footnote{Department of Mathematics, Sungkyunkwan University, Suwon, South Korea}, Sooyeong Kim\footnote{Department of Mathematics and Statistics, York University, Toronto, Canada.\label{foot1} (Contact: kimswim@yorku.ca)}, Minho Song\footnote{Applied Algebra and Optimization Research Center, Sungkyunkwan University, Suwon, South Korea} \par \bigskip
		
		\today
	\end{center}
\begin{abstract}
	On trees of fixed order, we show a direct relation between Kemeny's constant and Wiener index, and provide a new formula of Kemeny's constant from the relation with a combinatorial interpretation. Moreover, the relation simplifies proofs of several known results for extremal trees in terms of Kemeny's constant for random walks on trees. Finally, we provide various families of co-Kemeny's mates, which are two non-isomorphic connected graphs with the same Kemeny's constant, and we also give a necessary condition for a tree to attain maximum Kemeny's constant for trees with fixed diameter.
\end{abstract}



\section{Introduction and preliminaries}

Kemeny's constant for random walks on a graph has been not only exploited in many applications (\textit{e.g.} \cite{crisostomi2011google,obereigner2021markovian,yilmaz2020kemeny}), as the expected time for a random walker to travel between two randomly chosen vertices, but also studied recently in combinatorial settings to understand its insight  \cite{ciardo2020kemeny,faught20221,kim2022families,kirkland2016kemeny}. In this article, we focus on the combinatorial expression for Kemeny's constant in \cite{kirkland2016kemeny} (see \cite{kemeny1983finite} for the original definition of Kemeny's constant in the context of Markov chains) and study its relation to Wiener index, which has been widely studied in mathematical chemistry \cite{knor2015mathematical}. We point out a direct relation between them on trees, and we present several benefits from the relation.

A simple, undirected graph $G$ consists of a set of vertices $V(G)$ and a set of edges $E(G) \subseteq \{\{u, v\} \mid u, v \in V(G)\}$.
The \emph{order} of \( G \) is \( |V(G)| \). We use $m_G$ to denote the number of edges. 
Two vertices $u$ and $v$ are said to be \emph{adjacent} if there is an edge $\{u, v\}$ in $E(G)$; this is also denoted by $u\sim v$. 
The \emph{degree} of a vertex $v$, denoted $\deg_G(v)$, is the number of vertices adjacent to $v$ in $G$. 
For $v,w\in V(G)$, we use $\mathrm{dist}_G(v,w)$ to denote the distance between $v$ and $w$; that is, the length of the shortest path between $v$ and $w$. 
For $X\subseteq V(G)$, we define $\mathrm{dist}_G(X,v)=\min_{w\in X}\mathrm{dist}_G(w,v)$. 
The \emph{eccentricity} $e_G(v)$ of a vertex $v$ in $G$ is defined as $e_G(v)=\max_{w\in V(G)}\mathrm{dist}_G(w,v)$. 
The \emph{diameter} $\mathrm{diam}(G)$ of $G$ is defined as $\mathrm{diam}(G) = \max_{v\in V(G)}e(v)$. 
If the vertices in $V(G)$ are labelled $v_1, \ldots, v_n$, we define $\mathbf{d}_G$ to be the column vector whose $i^\text{th}$ component is $\mathrm{deg}_G(v_i)$ for $1\leq i\leq n$. 
We denote by $\mathbf{1}_k$ the all-ones vector of length $k$, by $\mathbf{0}_k$ the all-zeros vector of length $k$. 
The subscript $k$ is omitted if the size is clear from the context.

A \textit{tree} is a connected graph that has no cycles. A \textit{forest} is a graph whose connected components are trees. A \textit{spanning tree} (resp. a \textit{spanning forest}) of a graph $G$ is a subgraph that is a tree (resp. a forest) and includes all of the vertices of $G$. In particular, a \textit{$2$-tree spanning forest} of $G$ is a spanning forest that consists of $2$ trees.

Let $\cT$ be a tree. A vertex $c$ of $\cT$ is called a \emph{center} if $e(c) = \min_{v\in V(\cT)} e(v)$. We use $C(\cT)$ to denote the set of centers of $\cT$. 
Then, a center $c$ of \( \cT \) must be on any path of length $\mathrm{diam}(\cT)$. It follows that there is a unique center if $\mathrm{diam}(\cT)$ is even; and there are exactly two adjacent centers if $\mathrm{diam}(\cT)$ is odd. 

Let $G$ be a connected graph with \( V(G) = \{1,2,\dots,n\} \). We denote by $\tau_{G}$ the number of spanning trees of $G$, and by $f_{i,j}^{G}$ the number of $2$-tree spanning forests of $G$ such that one of the two trees contains vertex $i$ of $G$, and the other has vertex $j$ of $G$. We define $F_{G}$ to be the matrix given by $F_{G}=[f_{i,j}^{G}]_{i,j=1}^n$. Then, \emph{Kemeny's constant $\kappa(G)$ of a random walk on $G$} is given \cite{kirkland2016kemeny} by
\begin{equation}\label{formula:KemenyF}
\kappa(G)=\frac{\mathbf{d}_G^T F_G\mathbf{d}_G}{4m_G\tau_G}.
\end{equation}
We simply call $\kappa(G)$ \emph{Kemeny's constant of $G$}. Two non-isomorphic connected graphs $G_1$ and $G_2$ on the same number of vertices are said to be \textit{co-Kemeny mates} if $\kappa(G_1)=\kappa(G_2)$ (see \cite{kirkland2016kemeny}).

Let $G$ be a connected graph on \( V(G) = \{1,2,\dots,n\} \). The distance matrix $D$ of $G$ is given by $D=[d_{i,j}]_{i,j=1}^n$ where $d_{i,j}=\mathrm{dist}_G(i,j)$. The \textit{Gutman index} $\mathrm{Gut}(G)$ \cite{gutman1994selected} and the \textit{Wiener index} $W(G)$ \cite{wiener1947structural} of \( G \) are defined by 
\begin{align*}
\mathrm{Gut}(G)= \frac{1}{2}\mathbf{d}_G^T D\mathbf{d}_G, \quad\mbox{and}\quad W(G) = \frac{1}{2}\mathbf{1}^TD\mathbf{1}.
\end{align*}

\section{Relation between Kemeny's constant and Wiener index}

Let $\cT$ be a tree on the vertex set \( \{1,\dots,n\} \). Then, $F_\cT$ is the distance matrix of $\cT$, that is, $f_{i,j}^\cT = \mathrm{dist}_{\cT}(i,j)$. It is found in \cite[Lemma 8.7]{bapat2010graphs} that $\mathbf{d}_\cT^T F_\cT=\mathbf{1}^T(2F_\cT-(n-1)I)$. It follows that
\begin{align}\label{relation:kem and wiener}
\kappa(\cT) = \frac{4\mathbf{1}^TF_\cT\mathbf{1}-2(n-1)(2n-1)}{4(n-1)} = \frac{\mathbf{1}^TF_\cT\mathbf{1}}{n-1}-\frac{(2n-1)}{2} = \frac{2W(\cT)}{n-1}-n+\frac{1}{2}.
\end{align}
Therefore, ordering of two Kemeny's constants for two trees $\cT_1$ and $\cT_2$ on the same vertices is determined by comparing the sum of all entries in the distance matrix of $\cT_1$ with that of $\cT_2$. In other words, $\kappa(\cT_1)<\kappa(\cT_2)$ if and only if $W(\cT_1)<W(\cT_2)$; and $\kappa(\cT_1)=\kappa(\cT_2)$ if and only if $W(\cT_1)=W(\cT_2)$.

\begin{theorem}\label{thm:ordering}
	Let $\mathcal{C}$ be the set of trees on $n$ vertices satisfying some property. Then, trees with maximum/minimum Kemeny's constant in $\mathcal{C}$ are exactly the same as those with maximum/minimum Wiener index.
\end{theorem}

\begin{theorem}\label{prop:1}
	Two non-isomorphic trees $\cT_1$ and $\cT_2$ on the same number of vertices are co-Kemeny mates if and only if $W(\cT_1)=W(\cT_2)$.
\end{theorem}

We now provide an interpretation for Kemeny's constant on trees. Let $n_1(e)$ and $n_2(e)$ be the number of vertices in each of the two components of $\cT\backslash e$, respectively, where $\cT\backslash e$ is the forest obtained from $\cT$ by removing $e$. Then, we have a useful formula \cite{mohar1988compute}
\begin{align}\label{formula:sum over edges}
W(\cT) = \sum_{e\in E(\cT)} n_1(e)n_2(e).
\end{align}
Note $n_1(e)+n_2(e)=n$. Hence, one can find
\begin{align}\label{newformula}
\kappa(\cT) = \frac{1}{2}\frac{\sum_{e\in E(\cT)}(2n_1(e)-1)(2n_2(e)-1)}{n-1}.
\end{align}
Here the quantity $(2n_1(e)-1)(2n_2(e)-1)$ can be interpreted as follows. For $e=v_1\sim v_2$, consider two subtrees $\cT_1$ and $\cT_2$ in $\cT\backslash e$ containing $v_1$ and $v_2$, respectively. Form $\hat{\cT}$ from $\cT$ by attaching $v_i$ of a copy of $\cT_i$ to $v_i$ of $\cT$ for $i=1,2$. Then, the quantity is the number of shortest paths of $\hat{\cT}$ such that they pass through $e$. Therefore, $2\kappa(\cT)$ is the average of such quantities for all edges in $\cT$.

\begin{remark}
	The formula \eqref{newformula} reduces computation time of Kemeny's constant for trees. So, it can simplifies finding formulas of Kemeny's constant for families of particular trees. For instance, the formula for broom stars appears in \cite{ciardo2020kemeny} and one can deduce the same one by \eqref{newformula} with less computation. 
\end{remark}

In term of comparison of Kemeny's constants for trees on the same vertices, it is easier to compare their Wiener indices. Let us define a function $\omega$ from $E(\cT)$ to the set of positive numbers given by $\omega(e)=n_1(e)n_2(e)$. Then, $\omega(e)$ is the number of shortest paths of $\cT$ that pass through the specified edge $e\in E(\cT)$. We denote by $\cT^\omega$ the weighted tree obtained from $\cT$ in which each edge $e$ assigns $\omega(e)$. For two trees $\cT_1$ and $\cT_2$ on the same number of vertices, ordering of $\kappa(\cT_1)$ and $\kappa(\cT_2)$ is determined by the sums of all weights for the weighted trees $\cT_1^\omega$ and $\cT_2^\omega$.

\begin{example}
	Let us consider two trees $\cT_1$ and $\cT_2$ on $6$ vertices:
	
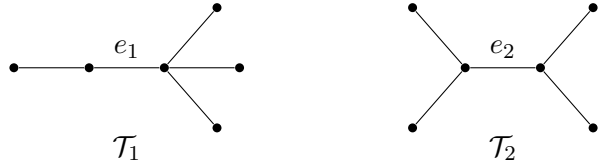
\begin{center}
		\begin{tikzpicture}
		\node[circle, fill, inner sep=1.2pt] (v1) at (-0.5,0) {};
		\node[circle, fill, inner sep=1.2pt] (v2) at (0.5,0) {};
		\node[circle, fill, inner sep=1.2pt] (v3) at (-1.5,0) {};
		\node[circle, fill, inner sep=1.2pt] (v4) at (1.2,0.8) {};
		\node[circle, fill, inner sep=1.2pt] (v5) at (1.5,0) {};
		\node[circle, fill, inner sep=1.2pt] (v6) at (1.2,-0.8) {};
		\node[label={$\cT_1$}] (v) at (0,-1.5) {};
		
		\draw (v1)--(v2) node [midway, above] {$e_1$};
		\draw (v3)--(v1);
		\draw (v2)--(v4);
		\draw (v2)--(v5);
		\draw (v2)--(v6);
		
		\begin{scope}[xshift=5cm]		
		\node[circle, fill, inner sep=1.2pt] (v1) at (-0.5,0) {};
		\node[circle, fill, inner sep=1.2pt] (v2) at (0.5,0) {};
		\node[circle, fill, inner sep=1.2pt] (v3) at (-1.2,0.8) {};
		\node[circle, fill, inner sep=1.2pt] (v4) at (-1.2,-0.8) {};
		\node[circle, fill, inner sep=1.2pt] (v5) at (1.2,0.8) {};
		\node[circle, fill, inner sep=1.2pt] (v6) at (1.2,-0.8) {};
		\node[label={$\cT_2$}] (v) at (0,-1.5) {};
		
		\draw (v1)--(v2) node [midway, above] {$e_2$};
		\draw (v3)--(v1);
		\draw (v6)--(v2)--(v5);
		\draw (v1)--(v4);
		\end{scope}
		\end{tikzpicture}
	\end{center}
	It would not be prompt to tell which one has bigger Kemeny's constant, considering \eqref{formula:KemenyF}, other formulae in the context of Markov chains, or a qualitative interpretation (how well-connected a graph is) in the context of random walks on undirected graphs. As discussed above, we only need to consider values of the function $\omega$ for all edges in $\cT_1$ and $\cT_2$.
	
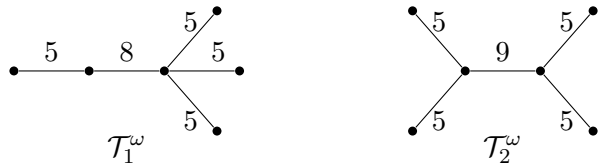
\begin{center}
		\begin{tikzpicture}
		\node[circle, fill, inner sep=1.2pt] (v1) at (-0.5,0) {};
		\node[circle, fill, inner sep=1.2pt] (v2) at (0.5,0) {};
		\node[circle, fill, inner sep=1.2pt] (v3) at (-1.5,0) {};
		\node[circle, fill, inner sep=1.2pt] (v4) at (1.2,0.8) {};
		\node[circle, fill, inner sep=1.2pt] (v5) at (1.5,0) {};
		\node[circle, fill, inner sep=1.2pt] (v6) at (1.2,-0.8) {};
		\node[label={$\cT_1^\omega$}] (v) at (0,-1.5) {};
		
		\draw (v1)--(v2) node [midway, above] {$8$};
		\draw (v1)--(v3) node [midway, above] {$5$};
		\draw (v2)--(v4) node [midway, above] {$5$};
		\draw (v2)--(v5) node [midway, above, near end] {$5$};
		\draw (v2)--(v6) node [midway, below] {$5$};
		
		\begin{scope}[xshift=5cm]		
		\node[circle, fill, inner sep=1.2pt] (v1) at (-0.5,0) {};
		\node[circle, fill, inner sep=1.2pt] (v2) at (0.5,0) {};
		\node[circle, fill, inner sep=1.2pt] (v3) at (-1.2,0.8) {};
		\node[circle, fill, inner sep=1.2pt] (v4) at (-1.2,-0.8) {};
		\node[circle, fill, inner sep=1.2pt] (v5) at (1.2,0.8) {};
		\node[circle, fill, inner sep=1.2pt] (v6) at (1.2,-0.8) {};
		\node[label={$\cT_2^\omega$}] (v) at (0,-1.5) {};
		
		\draw (v1)--(v2) node [midway, above] {$9$};
		\draw (v1)--(v3) node [midway, above] {$5$};
		\draw (v1)--(v4) node [midway, below] {$5$};
		\draw (v2)--(v5) node [midway, above] {$5$};
		\draw (v2)--(v6) node [midway, below] {$5$};
		\end{scope}
		\end{tikzpicture}
	\end{center}
	Note that any edge incident to a pendent vertex has weight $5$. Since $\omega(e_2)>\omega(e_1)$, we have $\kappa(\cT_2)>\kappa(\cT_1)$.
\end{example}

\begin{remark}
	The relation \eqref{relation:kem and wiener} does not hold for graphs with cycles in general. Consider the following two unicycles $U_1$ and $U_2$ with a cycle of the same length.
	
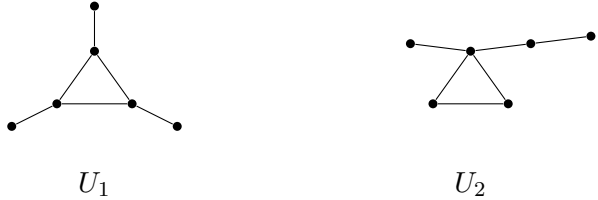
\begin{center}
		\begin{tikzpicture}
		\node[circle, fill, inner sep=1.2pt] (v1) at (-0.5,0) {};
		\node[circle, fill, inner sep=1.2pt] (v2) at (0.5,0) {};
		\node[circle, fill, inner sep=1.2pt] (v3) at (0,0.7) {};
		\node[circle, fill, inner sep=1.2pt] (v4) at (-1.1,-0.3) {};
		\node[circle, fill, inner sep=1.2pt] (v5) at (1.1,-0.3) {};
		\node[circle, fill, inner sep=1.2pt] (v6) at (0,1.3) {};
		\node[label={$U_1$}] (v) at (0,-1.5) {};
		
		\draw (v1)--(v2)--(v3)--(v1);
		\draw (v1)--(v4);
		\draw (v2)--(v5);
		\draw (v3)--(v6);
		
		\begin{scope}[xshift=5cm]		
		\node[circle, fill, inner sep=1.2pt] (v1) at (-0.5,0) {};
		\node[circle, fill, inner sep=1.2pt] (v2) at (0.5,0) {};
		\node[circle, fill, inner sep=1.2pt] (v3) at (0,0.7) {};
		\node[circle, fill, inner sep=1.2pt] (v4) at (-0.8,0.8) {};
		\node[circle, fill, inner sep=1.2pt] (v5) at (0.8,0.8) {};
		\node[circle, fill, inner sep=1.2pt] (v6) at (1.6,0.9) {};
		\node[label={$U_2$}] (v) at (0,-1.5) {};
		
		\draw (v1)--(v2)--(v3)--(v1);
		\draw (v3)--(v4);
		\draw (v3)--(v5);
		\draw (v5)--(v6);
		\end{scope}
		\end{tikzpicture}
	\end{center}
	They have the same Wiener index, but we have $\kappa(U_1)\approx 5.4167$ and $\kappa(U_2)\approx6.0833$ by \eqref{formula:KemenyF}. 
\end{remark}

\begin{remark}
	One may define Kemeny's constant for random walks on directed graphs (see \cite{kirkland2016kemeny} for the combinatorial expression) if they are strongly connected---that is, for each pair of vertices, there is a directed path between them. On the other hand, Wiener index can be defined for any directed graphs (see \cite{knor2016some}). We note that for trees, since any edge of a tree is a bridge, Kemeny's constant is only defined for undirected trees.
\end{remark}

Wiener index is closely related to Kemeny's constant for random walks on a tree. This relation has not attracted considerable attention yet, despite any extremal problems in terms of Kemeny's constants for trees on the same vertices is equivalent to those in terms of Wiener indices. So, we shall link known results in two different contexts.

Under several constraints, minimal Kemeny's constant on trees $\cT$ has been explored: \begin{enumerate}[label=(\alph*)]
	\item\label{case 1} for fixed order, the minimum of Kemeny's constant over all trees is uniquely attained at a star \cite{ciardo2020kemeny,kim2022families,kirkland2016kemeny} (they provide different proofs); and
	\item\label{case 2} for fixed order and diameter, the minimum of Kemeny's constant over all trees is uniquely attained at a particular caterpillar tree \cite{ciardo2020kemeny}.
\end{enumerate} 
For the proof of \ref{case 1}, \cite{ciardo2020kemeny} uses spectral property of graphs and the others \cite{kim2022families,kirkland2016kemeny} focus on minimizing $4m_\cT\kappa (\cT)$. The same result can be found in \cite{tomescu1999some} in the context of Gutman index. The proof of \ref{case 2} in \cite{ciardo2020kemeny} is also based on minimizing $4m_\cT\kappa (\cT)$.

Here we provide a different perspective. From \Cref{thm:ordering}, we may identify trees minimizing Kemeny's constant with trees minimizing Wiener index. Using \eqref{formula:sum over edges}, one can find that the minimum of Wiener index over trees on the same vertices is uniquely attained at a star, and hence \ref{case 1} follows. Furthermore, the authors of \cite{wang2008trees} present trees minimizing Wiener index for fixed order and diameter, using \eqref{formula:sum over edges} and some graph operations that result in decease of Wiener index. So, the result in \cite{wang2008trees} implies \ref{case 2}.

In contrast to the minimal Kemeny's constant or Wiener index problem, the maximal problem for trees with fixed order and diameter still remains an open problem. As a special case, trees attaining maximum of Wiener index over trees with fixed order and diameter at most $6$ are characterized in \cite{mukwembi2014wiener}.

Besides, minimal/maximal Wiener index problems under various constraints, which are equivalent to minimal/maximal Kemeny's constant problems, have been extensively studied. The authors of \cite{du2010minimum} show trees minimizing Wiener index for fixed order and matching number. We refer to the introduction of \cite{lin2014note} for citations of various results for trees with fixed order and maximum degree, for trees with all degrees odd, trees with fixed degree sequence, and so on. We also refer the interested reader to \cite{knor2015mathematical} for concrete open problems. 

\section{Some results in terms of Kemeny's constant from the relation}

In this section, we consider some operations on trees to manipulate Wiener index. Using them, we furnish families of co-Kemeny mates and give a necessary condition for trees with fixed order and diameter to attain maximum Kemeny's constant/Wiener index.

To find trees that minimize or maximize Wiener index under various conditions described in the previous section, several operations have been considered: generalized tree shift in \cite{scik2010on}; inner-moving, edge-growing, and lengthening transformation in \cite{wang2008trees}; $\alpha$--operation and $\beta$--operation in \cite{lin2014extremal}; diameter-growing transformation in \cite{lin2014note}; and relocating transformation in \cite{MR4292085}. 

Now, we consider two operations: one is a new operation, and the other appears in \cite{kirkland2016kemeny}, which is a generalized version of the operations in \cite{MR4292085,lin2014note}, inner-moving, lengthening transformation, and $\beta$--operation.

\begin{figure}
	\begin{center}
		\begin{tikzpicture}
		\tikzset{enclosed/.style={draw, circle, inner sep=0pt, minimum size=.10cm, fill=black}}
		\node[ label=\( \cT_1 \)] at (-2.8,-0.37) {};
		\node[enclosed] (-1) at (-1,0) {};
		\foreach \pos in {0,...,3} {
			\node[enclosed] (\pos) at (\pos,0) {} ;
		}
		\foreach \x in {0,...,2} {
			\draw plot [smooth cycle, tension=1] coordinates {(\x,0) (\x+0.12,0.2) (\x+0.4,1) (\x-0.4,1) (\x-0.12,0.2) };
		}
		\node[label=\( C_{1} \)] at (0,0.4) {};
		\node[label=\( \cdots \)] at (1,0.4) {};
		\node[label=\( C_{d-1} \)] at (2,0.4) {};
		\node[label=\( C_{0} \)] at (-1.7,-0.4) {};
		\node[label=\( C_{d} \)] at (3.7,-0.4) {};
		\draw plot [smooth cycle, tension=1] coordinates {(-1,0) (-1.2,0.1) (-2,0.2) (-2,-0.2) (-1.2,-0.1) };
		\draw plot [smooth cycle, tension=1] coordinates {(3,0) (3.2,0.1) (4,0.2) (4,-0.2) (3.2,-0.1) };
		\draw (-1) -- (3) ;
		
		\begin{scope}[shift={(0,-1.8)}]
		\node[ label=\( \cT_2 \)] at (-2.8,-0.37) {};
		\node[enclosed] (-1) at (-1,0) {};
		\foreach \pos in {0,...,3} {
			\node[enclosed] (\pos) at (\pos,0) {} ;
		}
		\foreach \x in {-1,...,1} {
			\draw plot [smooth cycle, tension=1] coordinates {(\x,0) (\x+0.13,0.2) (\x+0.4,1) (\x-0.4,1) (\x-0.1,0.2) };
		}
		\node[label=\( C_{1} \)] at (-1,0.4) {};
		\node[label=\( \cdots \)] at (0,0.4) {};
		\node[label=\( C_{d-1} \)] at (1,0.4) {};
		\node[label=\( C_{0} \)] at (-1.7,-0.4) {};
		\node[label=\( C_{d} \)] at (3.7,-0.4) {};
		\draw plot [smooth cycle, tension=1] coordinates {(-1,0) (-1.2,0.1) (-2,0.2) (-2,-0.2) (-1.2,-0.1) };
		\draw plot [smooth cycle, tension=1] coordinates {(3,0) (3.2,0.1) (4,0.2) (4,-0.2) (3.2,-0.1) };
		\draw (-1) -- (3) ;
		\end{scope}
		\end{tikzpicture}
		\caption{An illustration of the graph operation in \Cref{prop:operation1}.}
		\label{Figure:opertaion1}
	\end{center}
\end{figure}
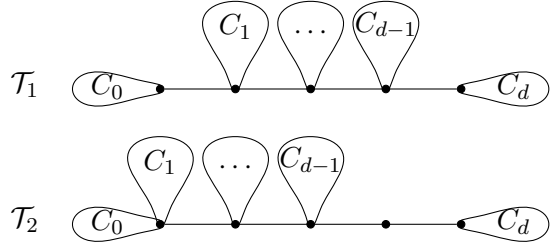

\begin{proposition}\label{prop:operation1}
	Let $\cT_1$ be a tree on \( n \) vertices with two distinct vertices $i_1$ and $i_2$. Denote the path of $\cT_1$ from $i_1$ to $i_2$ by $i_1\equiv l_0\sim l_1\sim \dots \sim l_d\equiv i_2$ for some $d\geq 2$. Form $\cT_2$ from $\cT_1$ by contracting the edge $\{l_0,l_1\}$ and subdividing the edge $\{l_{d-1},l_d\}$.
	
	Let $C_0$ (resp. \( C_d \)) denote the component of $\cT\backslash \{l_0\sim l_1\}$ (resp. $\cT\backslash \{l_{d-1}\sim l_d\}$) containing vertex $l_0$ (resp. \( l_d \)), and $C_i$ denote the component of $\cT\backslash \{l_{i-1}\sim l_i,l_{i}\sim l_{i+1}\}$ containing vertex $l_i$ for \( 1\le i \le d-1 \). (See Figure~\ref{Figure:opertaion1} for an illustration.) Then, 
	\begin{align*} 
	W(\cT_1)-W(\cT_2)&=|V(C_0)|(n-|V(C_0)|)-(|V(C_d)|+1)(n-|V(C_d)|-1) \\
	&+(d-2)(n+1)-2(d-2)|V(C_0)|-2\sum_{i=1}^{d-2}(d-1-i)|V(C_i)|.
	\end{align*}
\end{proposition}
\begin{proof}
	We use \eqref{formula:sum over edges} for the proof. 
	Let \( \omega_{i}(e) \) be the weight of edge \( e \) in \( \cT_i^\omega \) for \( i=1,2 \). Observe that the weight of each edge in \( C_i \) for \( 0\le i\le d \) is invariant under the operation described in Figure~\ref{Figure:opertaion1}. Thus it suffices to consider the weight of edges in the path from \( i_1 \) to \( i_2 \). Set \( \omega_i(j):=\omega_{i}(l_{j-1}\sim l_{j}) \) for brevity. Then,
	\[ 
	W(\cT_1)-W(\cT_2)=\sum_{j=1}^{d}\omega_1(j)-\sum_{j=1}^{d}\omega_2(j).
	\]
	Since \( \omega_1(d)=\omega_2(d) \) and \( \omega_1(j+1)-\omega_2(j)=n+1-2\sum_{i=0}^{j}|V(C_i)| \) for \( 1\le j\le d-2 \), we have
	\begin{align*} 
	&\sum_{j=1}^{d}\omega_1(j)-\sum_{j=1}^{d}\omega_2(j)\\
	&=\omega_1(1)-\omega_2(d-1)+\sum_{j=1}^{d-2}\left(n+1-2\sum_{i=0}^{j}|V(C_i)|\right)\\
	&=|V(C_0)|(n-|V(C_0)|)-(|V(C_d)|+1)(n-|V(C_d)|-1) \\
	&+(d-2)(n+1)-2(d-2)|V(C_0)|-2\sum_{i=1}^{d-2}(d-1-i)|V(C_i)|,
	\end{align*}
	which completes the proof.
\end{proof}

For some special cases, we obtain simple formulas.

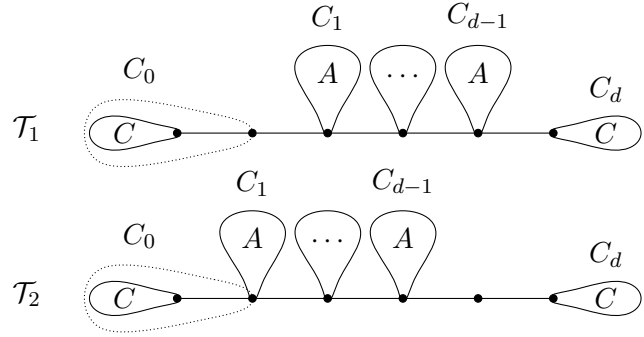
\begin{figure}
	\begin{center}
		\begin{tikzpicture}
		\tikzset{enclosed/.style={draw, circle, inner sep=0pt, minimum size=.10cm, fill=black}}
		\node[ label=\( \cT_1 \)] at (-4,-0.37) {};
		\node[enclosed] (-1) at (-1,0) {};
		\foreach \pos in {0,...,3} {
			\node[enclosed] (\pos) at (\pos,0) {} ;
		}
		\foreach \x in {0,...,2} {
			\draw plot [smooth cycle, tension=1] coordinates {(\x,0) (\x+0.12,0.2) (\x+0.4,1) (\x-0.4,1) (\x-0.12,0.2) };
		}
		\node[label=\( A \)] at (0,0.4) {};
		\node[label=\( C_1 \)] at (0,1.1) {};
		\node[label=\( \cdots \)] at (1,0.4) {};
		\node[label=\( A \)] at (2,0.4) {};
		\node[label=\( C_{d-1} \)] at (2,1.1) {};
		\node[label=\( C \)] at (-2.7,-0.4) {};
		\node[label=\( C \)] at (3.7,-0.4) {};
		\node[label=\( C_0 \)] at (-2.5,0.4) {};
		\node[label=\( C_d \)] at (3.7,0.15) {};
		\node[enclosed] at (-2,0) {};
		\draw[densely dotted] plot [smooth cycle, tension=0.75] coordinates {(-1,0) (-1.5,0.25) (-3,0.4) (-3,-0.4) (-1.5,-0.25) };
		\draw plot [smooth cycle, tension=1] coordinates {(-2,0) (-2.2,0.1) (-3,0.2) (-3,-0.2) (-2.2,-0.1) };
		\draw plot [smooth cycle, tension=1] coordinates {(3,0) (3.2,0.1) (4,0.2) (4,-0.2) (3.2,-0.1) };
		\draw (-1) -- (3) ;
		\draw (-1,0) -- (-2,0) ;
		
		\begin{scope}[shift={(0,-2.2)}]
		\node[ label=\( \cT_2 \)] at (-4,-0.37) {};
		\node[enclosed] (-1) at (-1,0) {};
		\foreach \pos in {0,...,3} {
			\node[enclosed] (\pos) at (\pos,0) {} ;
		}
		\foreach \x in {-1,...,1} {
			\draw plot [smooth cycle, tension=1] coordinates {(\x,0) (\x+0.13,0.2) (\x+0.4,1) (\x-0.4,1) (\x-0.1,0.2) };
		}
		\node[label=\( A \)] at (-1,0.4) {};
		\node[label=\( C_1 \)] at (-1,1.1) {};
		\node[label=\( \cdots \)] at (0,0.4) {};
		\node[label=\( A \)] at (1,0.4) {};
		\node[label=\( C_{d-1} \)] at (1,1.1) {};
		\node[label=\( C \)] at (-2.7,-0.4) {};
		\node[label=\( C \)] at (3.7,-0.4) {};
		\node[label=\( C_0 \)] at (-2.5,0.4) {};
		\node[label=\( C_d \)] at (3.7,0.15) {};
		\node[enclosed] at (-2,0) {};
		\draw[densely dotted] plot [smooth cycle, tension=0.75] coordinates {(-1,0) (-1.5,0.25) (-3,0.4) (-3,-0.4) (-1.5,-0.25) };
		\draw plot [smooth cycle, tension=1] coordinates {(-2,0) (-2.2,0.1) (-3,0.2) (-3,-0.2) (-2.2,-0.1) };
		\draw plot [smooth cycle, tension=1] coordinates {(3,0) (3.2,0.1) (4,0.2) (4,-0.2) (3.2,-0.1) };
		\draw (-1) -- (3) ;
		\draw (-1,0) -- (-2,0) ;
		\end{scope}
		\end{tikzpicture}
		\caption{An example of two isomorphic trees $\cT_1$ and $\cT_2$ for the operation in \Cref{Figure:opertaion1}, where $A$ and $C$ indicate some trees in their own oval boxes. Note \( |V(C_0)|-1=|V(C_d)|=|V(C)|\).}
		\label{Figure:iso}
	\end{center}
\end{figure}

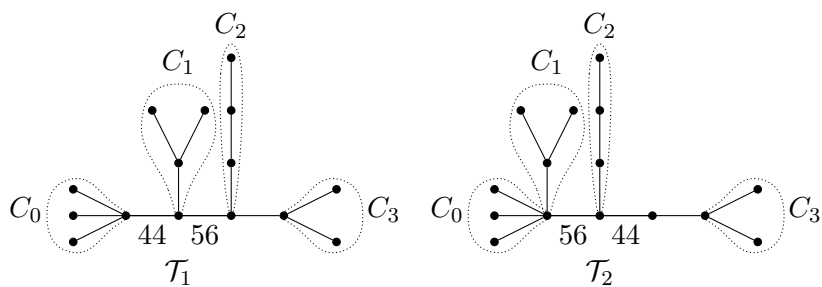
\begin{figure}
	\begin{center}
		\begin{tikzpicture}[scale = .7]
		\tikzset{enclosed/.style={draw, circle, inner sep=0pt, minimum size=.10cm, fill=black}}
		\node[label=\( \cT_1 \)] at (0,-1.7)	{};
		\node[label=\( C_0 \)] at (-2.9,-0.5)	{};
		\node[label=\( C_1 \)] at (0,2.3)	{};
		\node[label=\( C_2 \)] at (1,3)	{};
		\node[label=\( C_3 \)] at (3.9,-0.5)	{};

		\node[enclosed] (0) at (-1, 0) {};
		\node[enclosed] (1) at (0, 0) {};
		\node[enclosed] (2) at (1, 0) {};
		\node[enclosed] (12) at (1, 1) {};
		\node[enclosed] (13) at (1, 2) {};
		\node[enclosed] (14) at (1, 3) {};
		\node[enclosed] (3) at (2, 0) {};
		\node[enclosed] (9) at (3, 0.5) {};
		\node[enclosed] (11) at (3, -0.5) {};
		\node[enclosed] (4) at (-2, -0.5) {};
		\node[enclosed] (10) at (-2, 0) {};
		\node[enclosed] (5) at (-2, 0.5) {};
		\node[enclosed] (6) at (-0.5, 2) {};
		\node[enclosed] (7) at (0.5, 2) {};
		\node[enclosed] (8) at (0, 1) {};

		\draw (6) to (8);
		\draw (7) to (8);
		\draw (5) to (0);
		\draw (4) to (0);
		\draw (8) to (1);
		\draw (3) to (10);
		\draw (3) to (9);
		\draw (3) to (11);
		\draw (2) to (12);
		\draw (14) to (12);
		\draw (0) -- (1) node [midway, below] {$44$};
		\draw (1) -- (2) node [midway, below] {$56$};

		\draw[densely dotted] plot [smooth cycle, tension=1] coordinates {(-1,0) (-1.2,0.2) (-2,0.7) (-2.5,0) (-2,-0.7) (-1.2,-0.2)};
		\draw[densely dotted] plot [smooth cycle, tension=1] coordinates {(0,0) (-0.2,0.5) (-0.7,1.8) (0,2.5) (0.7,1.8) (0.2,0.5) };
		\draw[densely dotted] plot [smooth cycle, tension=1] coordinates {(2,0) (2.2,0.2) (3,0.7) (3.5,0) (3,-0.7) (2.2,-0.2) };
		\draw[densely dotted] plot [smooth cycle, tension=1] coordinates {(1,0) (0.8,2.7) (1.2,2.7)};
		
		\begin{scope}[shift={(8,0)}]
		\node[label=\( \cT_2 \)] at (0,-1.7)	{};
		\node[label=\( C_0 \)] at (-2.9,-0.5)	{};
		\node[label=\( C_2 \)] at (0,3)	{};
		\node[label=\( C_3 \)] at (3.9,-0.5)	{};

		\node[enclosed] (a0) at (-1, 0) {};
		\node[enclosed] (a1) at (0, 0) {};
		\node[enclosed] (a2) at (1, 0) {};
		\node[enclosed] (a3) at (2, 0) {};
		\node[enclosed] (a9) at (-2, 0) {};
		\node[enclosed] (a4) at (-2, -0.5) {};
		\node[enclosed] (a5) at (-2, 0.5) {};
		\node[enclosed] (a10) at (3, 0.5) {};
		\node[enclosed] (a11) at (3, -0.5) {};
		\node[enclosed] (a12) at (0, 1) {};
		\node[enclosed] (a13) at (0, 2) {};
		\node[enclosed] (a14) at (0, 3) {};
		
		\draw (a0) -- (a1) node [midway, below] {$56$};
		\draw (a1) -- (a2) node [midway, below] {$44$};

		\draw (a9) to (a3);
		\draw (a5) to (a0);
		\draw (a4) to (a0);
		\draw (a3) to (a10);
		\draw (a3) to (a11);
		\draw (a1) to (a12);
		\draw (a14) to (a12);
		
		\begin{scope}[shift={(-1,0)}]
		\node[label=\( C_1 \)] at (0,2.3)	{};
		
		\node[enclosed] (a6) at (-0.5, 2) {};
		\node[enclosed] (a7) at (0.5, 2) {};
		\node[enclosed] (a8) at (0, 1) {};
		
		\draw (a8) to (a0);
		\draw (a6) to (a8);
		\draw (a7) to (a8);
		
		\draw[densely dotted] plot [smooth cycle, tension=1] coordinates {(0,0) (-0.2,0.5) (-0.7,1.8) (0,2.5) (0.7,1.8) (0.2,0.5) };
		\draw[densely dotted] plot [smooth cycle, tension=1] coordinates {(1,0) (0.8,2.7) (1.2,2.7)};
		\end{scope}
		
		\draw[densely dotted] plot [smooth cycle, tension=1] coordinates {(-1,0) (-1.2,0.2) (-2,0.7) (-2.5,0) (-2,-0.7) (-1.2,-0.2)};
		\draw[densely dotted] plot [smooth cycle, tension=1] coordinates {(2,0) (2.2,0.2) (3,0.7) (3.5,0) (3,-0.7) (2.2,-0.2) };
		
		\end{scope}
		\end{tikzpicture}
		\caption{An example of co-Kemeny mates \( \cT_1 \) and \( \cT_2 \).}
		\label{fig:co-Kemeny op1}
	\end{center}	
\end{figure}

\begin{corollary}\label{cor:operation1}
	Under the same assumption as in Proposition~\ref*{prop:operation1}, let \( |V(C_1)|=\cdots=|V(C_{d-2})|=t \) and \( |V(C_0)|+m=|V(C_d)| \). Then, 
	\[ 
	W(\cT_1)-W(\cT_2)=d+(5-3d)t+2(d-2)|V(C_0)|-1-(t-1)(d-1)m.
	\]
	Furthermore, if \( |V(C_0)|=t \), then
	\begin{align}\label{eqn1}
	W(\cT_1)-W(\cT_2)=-(t-1)(d-1)(m+1).
	\end{align} 
\end{corollary}
If \( t=1 \), then we instantly have that \( \cT_1=\cT_2 \). We now consider the case \( t\ge2 \). From \eqref{eqn1}, we have two interesting findings.  If \( m=-1 \), then two trees $\cT_1$ and $\cT_2$ have the same Kemeny's constant so that if they are non-isomorphic, then they are co-Kemeny mates (see \Cref{Figure:iso} for isomorphic trees under the operation and see \Cref{fig:co-Kemeny op1} for an example of co-Kemeny's mates). Besides, if $m\geq 0$, then \( W(\cT_1)<W(\cT_2) \) since \( d,t\ge2 \). That is, Kemeny's constant increases under the operation corresponding to \eqref{eqn1}; so, this could be used to solve some extremal problems regarding Wiener index as done in \cite{wang2008trees}.

Now we consider a different operation, which is equivalent to the result \cite[Proposition 4.1]{kirkland2016kemeny} in terms of Kemeny's constants. With the operation, we provide a necessary condition for a tree with fixed order and diameter to have the maximum Kemeny's constant/Wiener index.

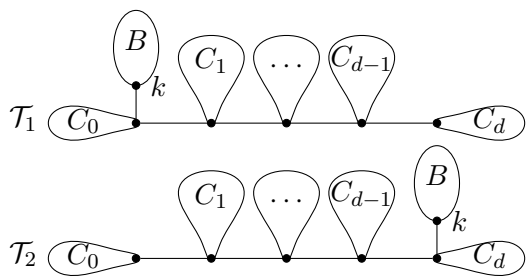
\begin{figure}
	\begin{center}
		\begin{tikzpicture}
		\tikzset{enclosed/.style={draw, circle, inner sep=0pt, minimum size=.10cm, fill=black}}
		\node[ label=\( \cT_1 \)] at (-2.5,-0.37) {};
		\node[enclosed] (-1) at (-1,0) {};
		\foreach \pos in {0,...,3} {
			\node[enclosed] (\pos) at (\pos,0) {} ;
		}
		\foreach \x in {0,...,2} {
			\draw plot [smooth cycle, tension=1] coordinates {(\x,0) (\x+0.12,0.2) (\x+0.4,1) (\x-0.4,1) (\x-0.12,0.2) };
		}
		\node[label=\( C_{1} \)] at (0,0.4) {};
		\node[label=\( \cdots \)] at (1,0.4) {};
		\node[label=\( C_{d-1} \)] at (2,0.4) {};
		\node[label=\( C_{0} \)] at (-1.7,-0.4) {};
		\node[label=\( C_{d} \)] at (3.7,-0.4) {};
		\draw plot [smooth cycle, tension=1] coordinates {(-1,0) (-1.2,0.1) (-2,0.2) (-2,-0.2) (-1.2,-0.1) };
		\draw plot [smooth cycle, tension=1] coordinates {(3,0) (3.2,0.1) (4,0.2) (4,-0.2) (3.2,-0.1) };
		\draw (-1) -- (3) ;
		
		\draw (-1,1) ellipse (0.3 and 0.5);
		\node[label = \( B \)] at (-1,0.7) {};
		\node[enclosed, label = {right: \( k \)}] at (-1,0.5) {};
		\draw (-1,0) -- (-1,0.5);
		
		\begin{scope}[shift={(0,-1.8)}]
		\node[ label=\( \cT_2 \)] at (-2.5,-0.37) {};
		\node[enclosed] (-1) at (-1,0) {};
		\foreach \pos in {0,...,3} {
			\node[enclosed] (\pos) at (\pos,0) {} ;
		}
		\foreach \x in {0,...,2} {
			\draw plot [smooth cycle, tension=1] coordinates {(\x,0) (\x+0.12,0.2) (\x+0.4,1) (\x-0.4,1) (\x-0.12,0.2) };
		}
		\node[label=\( C_{1} \)] at (0,0.4) {};
		\node[label=\( \cdots \)] at (1,0.4) {};
		\node[label=\( C_{d-1} \)] at (2,0.4) {};
		\node[label=\( C_{0} \)] at (-1.7,-0.4) {};
		\node[label=\( C_{d} \)] at (3.7,-0.4) {};
		\draw plot [smooth cycle, tension=1] coordinates {(-1,0) (-1.2,0.1) (-2,0.2) (-2,-0.2) (-1.2,-0.1) };
		\draw plot [smooth cycle, tension=1] coordinates {(3,0) (3.2,0.1) (4,0.2) (4,-0.2) (3.2,-0.1) };
		\draw (-1) -- (3) ;
		
		\draw (3,1) ellipse (0.3 and 0.5);
		\node[label = \( B \)] at (3,0.7) {};
		\node[enclosed, label = {right: \( k \)}] at (3,0.5) {};
		\draw (3,0) -- (3,0.5);
		\end{scope}
		\end{tikzpicture}
		\caption{An illustration of the graph operation in \Cref{prop:operation2}.}
		\label{Figure:opertaion2}
	\end{center}
\end{figure}

\begin{proposition}\cite{kirkland2016kemeny}\label{prop:operation2}
	Let $\cT$ be a tree with two distinct vertices $i_1$ and $i_2$, and $B$ be a tree with vertex $k$. For $j=1,2$, form $\cT_j$ from $\cT$ and $B$ by joining an edge between $i_j$ and $k$. 
	
	Denote the path of $\cT$ from $i_1$ to $i_2$ by $i_1\equiv l_0\sim l_1\sim \dots \sim l_d\equiv i_2$ for some $d\geq 1$. Let $C_0$ denote the component of $\cT\backslash \{l_0\sim l_1\}$ containing vertex $l_0$, and $C_d$ denote the component of $\cT\backslash \{l_{d-1}\sim l_d\}$ containing vertex $l_d$. When $d\geq 2$, we denote by $C_j$ the component of $\cT\backslash \{l_{j-1}\sim l_j,l_{j}\sim l_{j+1} \}$ containing vertex $l_j$ for $j=1,\dots,d-1$. (See Figure~\ref{Figure:opertaion2} for an illustration.) Then,
	\begin{align*}
	W(\cT_1)-W(\cT_2) = |V(B)|\sum_{j=0}^{d}|V(C_j)|(2j-d).
	\end{align*}
\end{proposition}

\begin{remark}
	Using the operation in \Cref{prop:operation2}, some families of co-Kemeny mates are provided in \cite{kirkland2016kemeny}.
\end{remark}

We shall consider a transformation from one tree to another, increasing Kemeny's constant while preserving its diameter.

\begin{definition}\label{def:cover}
	Let $\cT$ be a tree with two distinct vertices $i_1$ and $i_2$, and $B$ be a tree with vertex $k$. 
	Let $d = \mathrm{diam}(\cT)$.
	For $j=1,2$, form $\cT_j$ from $\cT$ and $B$ by joining an edge between $i_j$ and $k$.
	We say that \( \cT_1 \) \emph{covers} \( \cT_2 \), denoted \( \cT_2 \lessdot  \cT_1 \), if \( W(\cT_2) < W(\cT_1) \) and the diameters of \( \cT_1 \) and \( \cT_2 \) are the same.
	Let \( T_{n,d} \) be the set of trees on \( n \) vertices with diameter \( d \).
	We may consider \( T_{n,d} \) as a partially ordered set (or poset, for short) ordered by \( \cT \preceq \cT' \) in \( T_{n,d} \) if there is a sequence \( (\cT_0,\dots,\cT_k) \) of \( T_{n,d} \) such that \( \cT = \cT_0 \lessdot \dots \lessdot \cT_k = \cT' \).
\end{definition}

In \Cref{def:cover}, we can determine whether \( \cT_2 \lessdot \cT_1 \) or not by Proposition~\ref{prop:operation2}.
In particular, in the case of \( i_1 \sim i_2 \), maintaining the same notation in Proposition~\ref{prop:operation2}, we have
\begin{align}\label{eq:edge_translation}
\cT_1 \gtrdot \cT_2 \iff W(\cT_1)>W(\cT_2) \iff |V(C_0)|<|V(C_1)|.
\end{align}

We now give a necessary condition for maximal elements of the poset \( T_{n,d} \).

\begin{figure}
	\begin{center}
		\begin{tikzpicture}[scale = .5]
		\tikzset{enclosed/.style={draw, circle, inner sep=0pt, minimum size=.10cm, fill=black}}
		\node[enclosed, label={below : $u$}] (u) at (0,0) {};
		\node[enclosed, label={right : $k$}] (k) at (1,1) {};
		\node[enclosed, label={above : $v$}] (v) at (1,3.7) {};
		\node[enclosed, label={below : $i_1$}] (i_1) at (-4,-2) {};
		\node[enclosed, label={below : $i_2$}] (i_2) at (-2,-2) {};
		\node[enclosed, label={below : $i_j$}] (i_j) at (1,-2) {};
		\node[enclosed, label={below : $i_l$}] (i_l) at (4,-2) {};
		\node at (-0.5,-2) {\( \cdots \)};
		\node at (2.5,-2) {\( \cdots \)};
		\node at (0,-6) {\( \cT \)};
		
		\draw [bend right=15] (u) to (i_1);
		\draw [bend right=15] (u) to (i_2);
		\draw [bend left=15] (u) to (i_j);
		\draw [bend left=15] (u) to (i_l);
		\draw (u) to (k);
		
		\draw (i_1)+(0,-1.5) ellipse (0.8 and 1.5);
		\draw (i_2)+(0,-1.5) ellipse (0.8 and 1.5);
		\draw (i_j)+(0,-1.5) ellipse (0.8 and 1.5);
		\draw (i_l)+(0,-1.5) ellipse (0.8 and 1.5);
		\draw (k)+(0,1.5) ellipse (0.8 and 1.5);
		\node at (-4,-4) {\( C_1 \)};
		\node at (-2,-4) {\( C_2 \)};
		\node at (1,-4) {\( C_j \)};
		\node at (4,-4) {\( C_l \)};
		\node at (1,2.5) {\( \rotatebox{90}{\( \cdots \)} \)};
		
		\end{tikzpicture}\quad\qquad
		\begin{tikzpicture}[scale = .5]
		\tikzset{enclosed/.style={draw, circle, inner sep=0pt, minimum size=.10cm, fill=black}}
		\node[enclosed, label={below : $u$}] (u) at (0,0) {};
		\node[enclosed, label={right : $k$}] (k) at (2,0) {};
		\node[enclosed, label={above : $v$}] (v) at (2,2.7) {};
		\node[enclosed, label={below : $i_1$}] (i_1) at (-4,-2) {};
		\node[enclosed, label={below : $i_2$}] (i_2) at (-2,-2) {};
		\node[enclosed, label={below : $i_j$}] (i_j) at (1,-2) {};
		\node[enclosed, label={below : $i_l$}] (i_l) at (4,-2) {};
		\node at (-0.5,-2) {\( \cdots \)};
		\node at (2.5,-2) {\( \cdots \)};
		\node at (0,-6) {\( \cT_j \)};
		
		\draw [bend right=15] (u) to (i_1);
		\draw [bend right=15] (u) to (i_2);
		\draw [bend left=15] (u) to (i_j);
		\draw [bend left=15] (u) to (i_l);
		\draw [thick, bend right = 15] (i_j) to (k);
		
		\draw (i_1)+(0,-1.5) ellipse (0.8 and 1.5);
		\draw (i_2)+(0,-1.5) ellipse (0.8 and 1.5);
		\draw (i_j)+(0,-1.5) ellipse (0.8 and 1.5);
		\draw (i_l)+(0,-1.5) ellipse (0.8 and 1.5);
		\draw (k)+(0,1.5) ellipse (0.8 and 1.5);
		\node at (-4,-4) {\( C_1 \)};
		\node at (-2,-4) {\( C_2 \)};
		\node at (1,-4) {\( C_j \)};
		\node at (4,-4) {\( C_l \)};
		\node at (2,1.5) {\( \rotatebox{90}{\( \cdots \)} \)};
		\end{tikzpicture}
		\caption{An illustration for the proof of Proposition~\ref{thm:maximality}.}
		\label{Figure:Maximality}
	\end{center}
\end{figure}
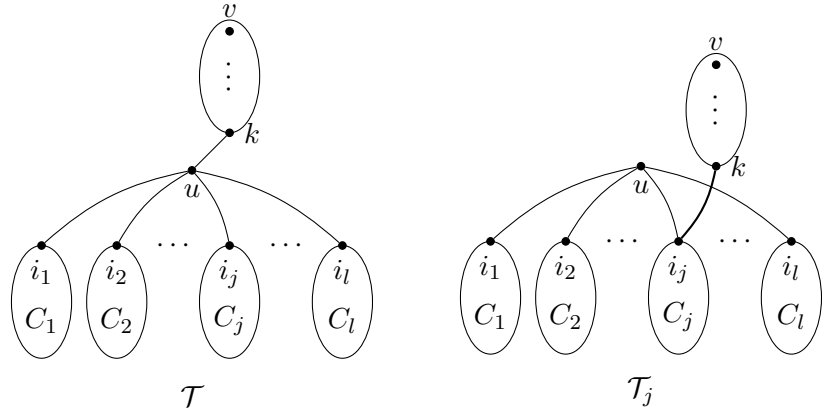

\begin{theorem}\label{thm:maximality}
	Let $n>d\geq 1$. If \( \cT\) is a maximal element of the poset \( T_{n,d} \), then 
	\[ \mathrm{dist}_{\cT}(C(\cT),v) =\lfloor d/2\rfloor \] 
	for all leaves \( v \) of \( \cT \).
\end{theorem}
\begin{proof}
	Assume to the contrary that there is a leaf \( v \) of \( \cT \) with \( \mathrm{dist}_{\cT}(C(\cT),v)<\lfloor d/2\rfloor \).
	Then, there exists a vertex of degree at least $3$. Choose the vertex $u$ nearest to $v$ such that \( \deg_\cT(u) = l+1 \) for some $l\geq 2$. Let \( \{i_1,\dots i_l,k\} \) be the set of neighbors of \( u \), where \( k \) is the nearest vertex to \( v \) among them (the vertex $k$ may be $v$). 
	Then, \( \dist_\cT(v,i_j)=\dist_\cT(v,k)+2 \) for \( j=1,\dots,l \).
	
	For each \( j=1,\dots,l \), we define \( C_j \) to be the component of \( \cT \backslash \{u\sim i_j\} \) containing \( i_j \), and also define \( \cT_j \) to be the tree obtained from $\cT$ by deleting $k\sim u$ and joining an edge between \( i_j \) and \( k \). (See Figure~\ref{Figure:Maximality} for an illustration.)
	Then, $\mathrm{dist}_{\cT_j}(C(\cT),v)-\mathrm{dist}_{\cT}(C(\cT),v)\leq 1$ and so $\mathrm{dist}_{\cT_j}(C(\cT),v)\leq \lfloor \frac{d}{2}\rfloor$. It follows that the diameters of \( \cT,\cT_1,\dots,\cT_l \) are the same, so \( \cT_1,\dots,\cT_l \in T_{n,d} \).
	
	For each \( j=1,\dots,l \), if \( W(\cT_j) > W(\cT) \), then it is a contradiction to the the assumption that \( \cT \) is a maximal element in $T_{n,d}$.
	Thus we have \( W(\cT_j)\le W(\cT) \), and by \eqref{eq:edge_translation},
	\begin{align*}
	|V(C_j)| 
	&\ge 1+ |V(C_1)|+\dots +|V(C_l)| - |V(C_j)| \\
	&> |V(C_1)|+\dots +|V(C_l)| - |V(C_j)|.
	\end{align*}
	Summing over all \( j \), we have 
	\begin{align*}
	\sum_{j=1}^l (|V(C_1)|+\dots +|V(C_l)| - 2|V(C_j)|) = (l-2)(|V(C_1)|+\dots +|V(C_l)|) < 0,
	\end{align*}
	which is impossible since \( l\ge 2 \). Therefore, our desired conclusion follows.
\end{proof}

By \Cref{thm:maximality}, to find the maximum Kemeny's constant on \( T_{n,d} \), it suffices to compare the Kemeny's constants of the maximal elements of \( T_{n,d} \).

\begin{example}\label{example: maximal}
	Let \( n=10 \) and \( d=4 \). There are 7 trees \( \cT_1,\dots,\cT_7 \) in \( T_{n,d} \) satisfying \( \mathrm{dist}_{\cT_i}(C(\cT_i),v) =\lfloor d/2\rfloor \) for all leaves \( v \) of \( \cT_i \) for each \( i=1,\dots,7 \) as follows:
	
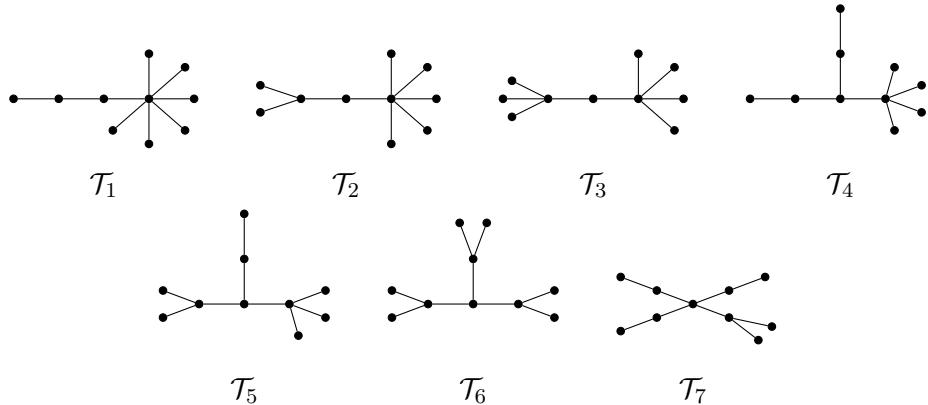
\begin{center}
		\begin{tikzpicture}[scale = 0.6]
		\tikzset{enclosed/.style={draw, circle, inner sep=1pt, minimum size=.10cm, fill=black}}
		\node[enclosed] (c) at (0,0) {};
		\node[enclosed] (i1) at (-1,0) {};
		\node[enclosed] (i2) at (1,0) {};
		\node[enclosed] (l1) at (-2,0) {};
		\node[enclosed] (l2) at (2,0) {};
		\node[enclosed] (l3) at (1.8,0.7) {};
		\node[enclosed] (l4) at (1,1) {};
		\node[enclosed] (l5) at (1.8,-0.7) {};
		\node[enclosed] (l6) at (1,-1) {};
		\node[enclosed] (l7) at (0.2,-0.7) {};
		\node[label={below, yshift=0cm: $\cT_1$}] (G) at (0,-1.2) {};
		
		\draw (c) -- (i1); 
		\draw (c) -- (i2); 
		\draw (i1) -- (l1);
		\draw (i2) -- (l2); 
		\draw (i2) -- (l3); 
		\draw (i2) -- (l4);
		\draw (i2) -- (l5); 
		\draw (i2) -- (l6); 
		\draw (i2) -- (l7); 
		\end{tikzpicture}\qquad
		\begin{tikzpicture}[scale = 0.6]
		\tikzset{enclosed/.style={draw, circle, inner sep=1pt, minimum size=.10cm, fill=black}}
		\node[enclosed] (c) at (0,0) {};
		\node[enclosed] (i1) at (-1,0) {};
		\node[enclosed] (i2) at (1,0) {};
		\node[enclosed] (l1) at (-1.9,0.3) {};
		\node[enclosed] (l2) at (2,0) {};
		\node[enclosed] (l3) at (1.8,0.7) {};
		\node[enclosed] (l4) at (1,1) {};
		\node[enclosed] (l5) at (1.8,-0.7) {};
		\node[enclosed] (l6) at (1,-1) {};
		\node[enclosed] (l7) at (-1.9,-0.3) {};
		\node[label={below, yshift=0cm: $\cT_2$}] (G) at (0,-1.2) {};
		
		\draw (c) -- (i1); 
		\draw (c) -- (i2); 
		\draw (i1) -- (l1);
		\draw (i2) -- (l2); 
		\draw (i2) -- (l3); 
		\draw (i2) -- (l4);
		\draw (i2) -- (l5); 
		\draw (i2) -- (l6); 
		\draw (i1) -- (l7); 
		\end{tikzpicture}\qquad
		\begin{tikzpicture}[scale = 0.6]
		\tikzset{enclosed/.style={draw, circle, inner sep=1pt, minimum size=.10cm, fill=black}}
		\node[enclosed] (c) at (0,0) {};
		\node[enclosed] (i1) at (-1,0) {};
		\node[enclosed] (i2) at (1,0) {};
		\node[enclosed] (l1) at (-2,0) {};
		\node[enclosed] (l2) at (2,0) {};
		\node[enclosed] (l3) at (1.8,0.7) {};
		\node[enclosed] (l4) at (1,1) {};
		\node[enclosed] (l5) at (-1.8,-0.4) {};
		\node[enclosed] (l6) at (1.8,-0.7) {};
		\node[enclosed] (l7) at (-1.8,0.4) {};
		\node[label={below, yshift=0cm: $\cT_3$}] (G) at (0,-1.2) {};
		
		\draw (c) -- (i1); 
		\draw (c) -- (i2); 
		\draw (i1) -- (l1);
		\draw (i2) -- (l2); 
		\draw (i2) -- (l3); 
		\draw (i2) -- (l4);
		\draw (i1) -- (l5); 
		\draw (i2) -- (l6); 
		\draw (i1) -- (l7); 
		\end{tikzpicture}\qquad
		\begin{tikzpicture}[scale = 0.6]
		\tikzset{enclosed/.style={draw, circle, inner sep=1pt, minimum size=.10cm, fill=black}}
		\node[enclosed] (c) at (0,0) {};
		\node[enclosed] (i1) at (0,1) {};
		\node[enclosed] (i2) at (-1,0) {};
		\node[enclosed] (i3) at (1,0) {};
		\node[enclosed] (l1) at (0,2) {};
		\node[enclosed] (l2) at (-2,0) {};
		\node[enclosed] (l3) at (1.2,0.7) {};
		\node[enclosed] (l4) at (1.8,0.3) {};
		\node[enclosed] (l5) at (1.8,-0.3) {};
		\node[enclosed] (l6) at (1.2,-0.7) {};
		\node[label={below, yshift=0cm: $\cT_4$}] (G) at (0,-1.2) {};
		
		\draw (c) -- (i1); 
		\draw (c) -- (i2); 
		\draw (c) -- (i3);
		\draw (i1) -- (l1); 
		\draw (i2) -- (l2); 
		\draw (i3) -- (l3);
		\draw (i3) -- (l4); 
		\draw (i3) -- (l5); 
		\draw (i3) -- (l6); 
		\end{tikzpicture}\qquad
		\begin{tikzpicture}[scale = 0.6]
		\tikzset{enclosed/.style={draw, circle, inner sep=1pt, minimum size=.10cm, fill=black}}
		\node[enclosed] (c) at (0,0) {};
		\node[enclosed] (i1) at (0,1) {};
		\node[enclosed] (i2) at (-1,0) {};
		\node[enclosed] (i3) at (1,0) {};
		\node[enclosed] (l1) at (0,2) {};
		\node[enclosed] (l2) at (-1.8,-0.3) {};
		\node[enclosed] (l3) at (-1.8,0.3) {};
		\node[enclosed] (l4) at (1.8,0.3) {};
		\node[enclosed] (l5) at (1.8,-0.3) {};
		\node[enclosed] (l6) at (1.2,-0.7) {};
		\node[label={below, yshift=0cm: $\cT_5$}] (G) at (0,-1.2) {};
		
		\draw (c) -- (i1); 
		\draw (c) -- (i2); 
		\draw (c) -- (i3);
		\draw (i1) -- (l1); 
		\draw (i2) -- (l2); 
		\draw (i2) -- (l3);
		\draw (i3) -- (l4); 
		\draw (i3) -- (l5); 
		\draw (i3) -- (l6); 
		\end{tikzpicture}\qquad
		\begin{tikzpicture}[scale = 0.6]
		\tikzset{enclosed/.style={draw, circle, inner sep=1pt, minimum size=.10cm, fill=black}}
		\node[enclosed] (c) at (0,0) {};
		\node[enclosed] (i1) at (0,1) {};
		\node[enclosed] (i2) at (-1,0) {};
		\node[enclosed] (i3) at (1,0) {};
		\node[enclosed] (l1) at (-0.3,1.8) {};
		\node[enclosed] (l2) at (-1.8,-0.3) {};
		\node[enclosed] (l3) at (-1.8,0.3) {};
		\node[enclosed] (l4) at (1.8,0.3) {};
		\node[enclosed] (l5) at (1.8,-0.3) {};
		\node[enclosed] (l6) at (0.3,1.8) {};
		\node[label={below, yshift=0cm: $\cT_6$}] (G) at (0,-1.2) {};
		
		\draw (c) -- (i1); 
		\draw (c) -- (i2); 
		\draw (c) -- (i3);
		\draw (i1) -- (l1); 
		\draw (i2) -- (l2); 
		\draw (i2) -- (l3);
		\draw (i3) -- (l4); 
		\draw (i3) -- (l5); 
		\draw (i1) -- (l6); 
		\end{tikzpicture}\qquad
		\begin{tikzpicture}[scale = 0.6]
		\tikzset{enclosed/.style={draw, circle, inner sep=1pt, minimum size=.10cm, fill=black}}
		\node[enclosed] (c) at (0,0) {};
		\node[enclosed] (i1) at (0.8,0.3) {};
		\node[enclosed] (i2) at (0.8,-0.3) {};
		\node[enclosed] (i3) at (-0.8,0.3) {};
		\node[enclosed] (i4) at (-0.8,-0.3) {};
		\node[enclosed] (l1) at (1.6,0.6) {};
		\node[enclosed] (l2) at (1.75,-0.5) {};
		\node[enclosed] (l3) at (1.45,-0.8) {};
		\node[enclosed] (l4) at (-1.6,0.6) {};
		\node[enclosed] (l5) at (-1.6,-0.6) {};
		\node[label={below, yshift=0cm: $\cT_7$}] (G) at (0,-1.2) {};
		
		\draw (c) -- (i1); 
		\draw (c) -- (i2); 
		\draw (c) -- (i3);
		\draw (c) -- (i4); 
		\draw (i1) -- (l1); 
		\draw (i2) -- (l2);
		\draw (i2) -- (l3); 
		\draw (i3) -- (l4); 
		\draw (i4) -- (l5); 
		\end{tikzpicture}
	\end{center}
	Among them, \( \cT_3 \), \( \cT_6 \), and \( \cT_7 \) are the maximal elements of \( T_{10,4} \), and the Wiener indices of \( \cT_3 \), \( \cT_6 \), \( \cT_7 \) are \( 112,117,114 \), respectively.
	Therefore, the maximum Kemeny's constant in \( T_{10,4} \) is attained at the tree \( \cT_6 \).
\end{example}

On trees with fixed order $n$ and diameter $d$, trees attaining the maximum Kemeny's constant (equivalently, Wiener index) have not been completely characterized so far; in the context of Wiener index, it is a long-standing open problem. As seen in \Cref{example: maximal}, we are able to reduce much time for comparisons of Kemeny's constants, using Wiener index. Moreover, as \Cref{thm:maximality} provides a characterization of maximal elements of $T_{n,d}$, we can narrow our focus to particular trees to find our desired trees.


\end{document}